\documentclass{amsart}
\usepackage{graphicx,enumerate} 
\usepackage{xcolor}
\usepackage{amssymb,latexsym,amsmath,amsthm,listings}
\usepackage{hyperref}
\usepackage{natbib}

\hypersetup{
    colorlinks=true,
    citecolor=blue,
    linkcolor=blue,
    urlcolor=blue
}
\usepackage[a4paper, margin=1in]{geometry}

\makeatletter
\renewcommand\section{\@startsection {section}{1}{\z@}
{-30pt \@plus -1ex \@minus -.2ex}
{2.3ex \@plus.2ex}
{\normalfont\normalsize\bfseries\boldmath}}
\renewcommand\subsection{\@startsection{subsection}{2}{\z@}
{-3.25ex\@plus -1ex \@minus -.2ex}
{1.5ex \@plus .2ex}
{\normalfont\normalsize\bfseries\boldmath}}
\renewcommand{\@seccntformat}[1]{\csname the#1\endcsname. }

\makeatother

\title[$1729$: The only Carmichael number of the form $2^np^3+1$]{One more curious property of $1729$: The only Carmichael number of the form $2^np^3+1$ with prime $p$}
\author{Jean Lucien Randrianantenaina and Ali K. Uncu}
\address{University of Bath, Faculty of Science, Department of Computer Science, Claverton Down, BA2 7AY, United Kingdom}
\email{\{jlr85,\ aku21\}@bath.ac.uk}
\date{\today}

\dedicatory{In memory of Florian Luca, a great inspiration and teacher.}

\newtheorem{theorem}{Theorem}
\newtheorem{lemma}{Lemma}
\newtheorem{proposition}{Proposition}
\newtheorem{corollary}{Corollary}
\theoremstyle{definition}

\newcommand{\K}{\mathcal{K}}
\newcommand{\ord}{\operatorname{ord}}


\begin{document}

\begin{abstract}
Florian Luca recently proved that if $2^np^3+1$ is a Carmichael number for some prime $p$ and integer $n\ge1$, then $p=3$. Here we complete the cubic classification by proving that $1729$, Hardy--Ramanujan taxicab number, is the only Carmichael number of the form $2^np^3+1$ with $p$ prime and $n\ge1$. The proof identifies the unique possible multiplicatively dependent prime factor, sharpens the known bound on Fermat prime factors of $2^nk+1$ from $k^2$ to $2k$, derives the explicit bound $n\le2742$, and concludes with a short sieve requiring neither factorisation nor primality testing.
\end{abstract}

\maketitle
\section{Introduction}

A Carmichael number is a composite integer $N$ such that $a^N\equiv a\pmod N$ for every integer $a$. Korselt's criterion \cite{Korselt} gives the following characterisation: a composite integer $N$ is Carmichael if and only if $N$ is square-free and $q-1\bigl\vert  N-1$ for every prime factor $q$ of $N$. In particular, every Carmichael number is odd and has at least three distinct prime factors.

Cilleruelo, Luca and Pizarro-Madariaga \cite{CLP} proved that for every fixed odd $k$, the sequence $(2^nk+1)_{n\ge1}$ contains only finitely many Carmichael numbers, and that the smallest element of the set
$$
\K:=\{k\text{ odd}:2^nk+1\text{ is Carmichael for some }n\ge1\}
$$
is $27$. In particular, for $n=6$ and $k=27$, we have $1729=2^6\,27+1=7\,13\,19$.

Alahmadi and Luca \cite{AL} subsequently proved that $\K$ contains no prime, and Luca and Randrianantenaina \cite{LR} proved that $\K$ contains no square of a prime. The cubic case is different, since $3^3=27$ belongs to $\K$.

Luca \cite{Luca1729} addressed the prime part of the cubic case and proved the following result.

\begin{theorem}[Luca \cite{Luca1729}]\label{thm:prime-part}
If $p$ is prime and $2^np^3+1$ is Carmichael for some $n\ge1$, then $p=3$.
\end{theorem}

Thus any Carmichael number of the form $2^np^3+1$ must occur in the sequence $2^n3^3+1=27\,2^n+1$. Our main result completes the cubic classification.

\begin{theorem}\label{thm:taxicab}
For $n\ge1$, $2^n3^3+1$ is a Carmichael number if and only if $n=6$. Consequently,
$$
1729=2^6\,3^3+1
$$
is the unique Carmichael number in the sequence $(2^n3^3+1)_{n\ge1}$.
\end{theorem}

Combining \autoref{thm:prime-part} and \autoref{thm:taxicab} gives the complete classification.

\begin{corollary}\label{cor:cubic}
For a prime $p$ and integer $n\ge1$, $2^np^3+1$ is Carmichael if and only if $(p,n)=(3,6)$.
\end{corollary}

In 2026, Luca \cite{Luca2026} proved that, for every fixed odd $m\ge5$, there are only finitely many primes $p$ with $p^m\in\K$; the resulting bound on $p$ is ineffective. Tiebekabe \cite{Tiebekabe} made several steps of \cite{CLP,Luca2026} explicit. In particular, if $m\ge3$ is odd, $p>m$ and $2^np^m+1$ is Carmichael, then $n<25m^4\log p$, and every prime factor $q$ of $2^np^m+1$ has $q-1$ multiplicatively independent of $2^np^m$ \cite[Theorems~1.1 and~1.3]{Tiebekabe}. He also reported an exhaustive search showing that $1729$ is the only Carmichael number $2^np^m+1\le10^{18}$ with $m$ odd and $3\le m\le37$ \cite[Theorem~1.4]{Tiebekabe}.

Both structural results just quoted require $p>m$, so they do not apply to the case $p=m=3$ treated here. This exclusion is genuine: for $1729$, the prime factor $13$ satisfies $13-1=12$ and $1729-1=12^3$, so a multiplicatively dependent prime factor occurs; see \autoref{lem:dependent}. Moreover, the restriction $27\,2^n+1\le10^{18}$ covers only $n\le55$, whereas \autoref{thm:taxicab} covers all $n\geq 1$.

The organization of this paper is as follows. In \autoref{sec2}, we define and classify the dependent and independent prime factors of $N_n$. \autoref{sec3} presents a sharper version of the bound in \cite[Lemma~2]{CLP} for Fermat prime factors. In \autoref{sec:bound}, we first bound the largest exponent $n$ that can yield a Carmichael number $N_n$, and then we give explicit bounds for the size of prime factors of $N_n$. \autoref{sec:sieve} has the sieve for the divisors of $N_n$ and this is where we also present the proofs of \autoref{thm:taxicab} and \autoref{cor:cubic}. Finally, we include an outlook section where we discuss possible future directions.

\section{Prime Factors and Dependency}\label{sec2}

Let $N_n:=2^n3^3+1=27\,2^n+1$. Assume from now on, unless stated otherwise, that $N_n$ is Carmichael. If $q\bigl\vert  N_n$ is prime, Korselt's criterion gives $q-1\bigl\vert  N_n-1=2^n3^3$. Therefore every prime factor of $N_n$ has the form
\begin{equation}\label{eq:factor}
q=2^a3^b+1,\qquad 1\le a\le n,\qquad 0\le b\le3.
\end{equation}

We call the integers $q-1$ and $N_n-1$ \textit{multiplicatively dependent}, if there exist positive integers $u,v$ such that
\begin{equation}\label{eq:dependence}
(2^a3^b)^u=(2^n3^3)^v.
\end{equation}

For $b=0$ we have $q=2^a+1$, so $a$ must be a power of $2$ and $q$ is a Fermat prime. For $b\ge1$, we now determine when $q-1=2^a3^b$ and $N_n-1=2^n3^3$ can be multiplicatively dependent. In the analogous cubic problem with $p>3$, Luca excludes this possibility \cite[Section~2.1]{Luca1729}; for $p=3$ it does occur, as the factor $13$ of $1729$ shows.

Here, we recite a known result for our cases \cite[Lemma 3]{CLP}.

\begin{lemma}\label{lem:dependent}
Let $q=2^a3^b+1\bigl\vert  N_n$ be prime, with $1\le b\le3$. If $N_n-1$ and $q-1$ are multiplicatively dependent, then $3\bigl\vert  n$, $b=1$ and $a=n/3$. Thus the only possible dependent prime factor is
\begin{equation}\label{eq:dependent-factor}
D_n:=3\,2^{n/3}+1.
\end{equation}
Conversely, if $3\bigl\vert  n$, then $D_n\bigl\vert  N_n$; if $D_n$ is prime, it is the unique dependent prime factor of $N_n$.
\end{lemma}

\begin{proof}
If $q-1$ and $N_n-1$ are multiplicatively dependent, comparing the $2$-adic and $3$-adic valuations in \eqref{eq:dependence} gives $au=nv$ and $bu=3v$, hence
\begin{equation}\label{eq:dependence-condition}
3a=nb.
\end{equation}
Conversely, \eqref{eq:dependence-condition} gives $(2^a3^b)^3=(2^n3^3)^b$, so it is equivalent to multiplicative dependence in the present setting.

If $b=3$, then \eqref{eq:dependence-condition} gives $a=n$ and $q=N_n$. This is a contradiction since $q$ is prime whereas $N_n$ is Carmichael and therefore composite.

If $b=2$, then $3a=2n$, so $n=3t$ and $a=2t$ for some $t\ge1$. Put $x:=3\,2^t$. Then $q=x^2+1$ and $N_n=x^3+1$. Since $q\bigl\vert  N_n$ and $x^2\equiv-1\pmod q$, we have $x^3+1\equiv-x+1\pmod q$. Thus $q\bigl\vert  x-1$, impossible because $0<x-1<x^2+1=q$.

Finally, if $b=1$, then \eqref{eq:dependence-condition} gives $n=3a$. If $n=3t$ and $x:=3\,2^t$, then
$$
N_n=x^3+1=(x+1)(x^2-x+1),
$$
so $D_n=x+1=3\,2^t+1$ divides $N_n$. If $D_n$ is prime, then $D_n-1=3\,2^t$ and $N_n-1=(3\,2^t)^3=(D_n-1)^3$, proving dependence. Uniqueness follows because $D_n$ is completely determined by $n$.
\end{proof}

Note that $b=0$ possibility for the prime divisors of $N_n$ is ignored in \autoref{lem:dependent}. This is because $b=0$ forces $v=0$, and later $a u=0$. This does not fit into the multiplicative dependence definition with positive $u$ and $v$. We handle $b=0$ case next. 

\section{Eliminating Fermat Prime Factors}\label{sec3}

By \cite[Lemma~2]{CLP}, every Fermat prime factor $q$ of $2^nk+1$ is less than $k^2$. For $k=27$ this gives the preliminary bound $q<729$. The following elementary theorem sharpens $k^2$ to $2k$, and yields $q<54$ in the present case.

\begin{theorem}\label{lem:fermat-bound}
Let $k\ge3$ be odd, $n\ge0$, and let $q=2^{2^\alpha}+1$ be a prime dividing $2^nk+1$. Then $q<2k$.
\end{theorem}

\begin{proof}
If $\alpha=0$, then $q=3<6\le2k$. Assume $\alpha\ge1$, put $c:=2^\alpha$, so that $2^c\equiv-1\pmod q$, and write $n=ct+r$ with $0\le r<c$. From $q\bigl\vert 2^nk+1$ we obtain $0\equiv(-1)^t2^rk+1\pmod q$. Multiplying by $2^{c-r}$ and using $2^c\equiv-1\pmod q$ gives
\begin{equation}\label{eq:k-congruence}
k\equiv(-1)^t2^{c-r}\pmod q.
\end{equation}
Assume, for a contradiction, that $q\ge2k$. Then $0<k<q$ and $2\le2^{c-r}\le2^c<q$.

If $t$ is even, \eqref{eq:k-congruence} forces $k\geq 3$ to be even and yields a contradiction. If $t$ is odd, \eqref{eq:k-congruence} forces $k=q-2^{c-r}$. If $r=0$, then $k=q-2^c=1$, a contradiction. If $r\ge1$, then $2^{c-r}\le2^{c-1}=(q-1)/2$, so $k\ge(q+1)/2>q/2$, contradicting $q\ge2k$.
\end{proof}

Thus, for $k=27$, the only possible Fermat prime factors of $N_n$ are $3$, $5$ and $17$.

\autoref{lem:fermat-bound} has several useful consequences. First off, the bound $q<2k$ is attained up to one unit for every Fermat prime $q\ge5$. Indeed, let
$q=2^c+1$
be a Fermat prime with $c\ge2$, and set $k:=(q+1)/2=2^{c-1}+1$ and $n:=c+1$. Then $k$ is odd and
$$
2^nk+1=2^{c+1}(2^{c-1}+1)+1=(2^c+1)^2=q^2.
$$
Thus $q\bigl\vert 2^nk+1$ and $q=2k-1$.

Moreover, the theorem also improves the bound for the entire Fermat part. Let $F_j:=2^{2^j}+1$ be the Fermat numbers, and suppose that $F_j$ is the largest Fermat prime dividing $2^nk+1$. Since $F_0F_1\,s F_{j-1}=F_j-2$, we have $F_0F_1\,s F_j=(F_j-1)^2-1<(F_j-1)^2$. Therefore, the product $F$ of the distinct Fermat prime divisors of $2^nk+1$ satisfies
$$
F<4k^2.
$$
For comparison, \cite[(15)]{CLP} obtains $F<k^4$ from the earlier individual bound $F_j<k^2$. In Luca's cubic setting $k=p^3$, this appears as $F<p^{12}=k^4$; see \cite[Section~2.1]{Luca1729}. Thus \autoref{lem:fermat-bound} improves the corresponding general Fermat part estimate from $k^4$ to $4k^2$.

To exemplify this, we look at our case $k=27$. The earlier bound $q<729$ leaves $257$ as an additional possible Fermat prime factor that must be excluded separately. Indeed, $27^{-1}\equiv-19\pmod{257}$, so $257\bigl\vert 27\,2^n+1$ would imply $2^n\equiv19\pmod{257}$, but $\ord_{257}(2)=16$, where $\ord_r(s)$ is the smallest positive integer $u$ such that $s^u \equiv 1 \pmod{r}$. A direct check of the $j=1,2,\dots,16$ powers $2^j$ modulo $257$ shows that $19\notin\langle2\rangle_{257}$, where $\langle2\rangle_{257}$ denote the orbit of powers of 2 modulo 257. Hence, $257\nmid N_n$. The sharper bound $q<54$ makes this additional check unnecessary.

Now we further investigate some properties of $N_n$.
\begin{lemma}\label{lem:fermat-divisibility}
Let $N_n=27\,2^n+1$. Then 
\begin{enumerate}[i.] \item $3\nmid N_n$; \item $5\bigl\vert  N_n$ if and only if $n\equiv1\pmod4$; and \item $17\nmid N_n$.
\end{enumerate}
\end{lemma}

\begin{proof}
\begin{enumerate}[i.]
    \item Obvious since $3\bigl\vert 27$.
    \item We have $5\bigl\vert  N_n$ if and only if $2^{n+1}\equiv-1\equiv2^2\pmod5$, since $27\equiv2\pmod5$. As $\ord_5(2)=4$, this is equivalent to $n\equiv1\pmod4$.
    \item Since $27\equiv10\pmod{17}$ and $10^{-1}\equiv12 \equiv -5\pmod{17}$, we have $17\bigl\vert  N_n$ if and only if $2^n\equiv5\pmod{17}$. However, $\ord_{17}(2)=8$ and the orbit of powers of 2 modulo 17 is $\langle2\rangle_{17}=(1\ 2 \ 4\ 8\ 16\ 15\ 13\ 9)\pmod{17}$, which does not contain $5$.
\end{enumerate} 
\end{proof}

\begin{theorem}\label{thm:no-fermat}
If $N_n$ is Carmichael, then $N_n$ has no Fermat prime factor. In particular, every prime divisor of $N_n$ has $b\ge1$ in \eqref{eq:factor}, and $n\not\equiv1\pmod4$.
\end{theorem}

\begin{proof}
By \autoref{lem:fermat-bound} and \autoref{lem:fermat-divisibility}, the only Fermat prime that can divide $N_n$ is $5$. Indeed, $3\nmid N_n$ and $17\nmid N_n$. Suppose that $5\bigl\vert  N_n$. Since $N_n$ is square-free, $5$ occurs exactly once in its prime factorisation. Therefore every other prime divisor $q$ of $N_n$ is not a Fermat prime and, by \eqref{eq:factor}, has the form $q=2^a3^b+1$ with $b\ge1$, implying $q\equiv1\pmod3$.
Reducing the prime factorisation of $N_n$ modulo $3$ therefore gives $N_n\equiv5\equiv2\pmod3$, contradicting $N_n\equiv1\pmod3$. Hence, $N_n$ has no Fermat prime factor. The final assertion follows from part~(ii) of \autoref{lem:fermat-divisibility}.
\end{proof}

\section{Bounding \texorpdfstring{$n$}{n}}\label{sec:bound}

Here we mostly follow the footsteps of \cite{CLP} to find a sharper bound for $n$ of $N_n$ and the exponents of 2 that appear in the representation of the prime divisors of $N_n$. This enables us to limit the number of cases to consider for a proof of \autoref{thm:taxicab}.

Let $R$ denote the product of all prime factors $q=2^a3^b+1$ of $N_n$ such that $q-1$ and $N_n-1$ are multiplicatively independent. By \autoref{thm:no-fermat}, all these factors have $b\ge1$. By \autoref{lem:dependent}, at most one prime factor, namely $D_n$, is omitted from $R$, and this can happen only when $3\bigl\vert  n$. Hence,
\begin{equation}\label{eq:R-versus-N}
R=N_n\quad\text{if }3\nmid n,\qquad\text{and}\qquad R\ge\frac{N_n}{D_n}\quad\text{if }3\bigl\vert  n.
\end{equation}

\begin{lemma}\label{lem:order-relation}
Let $q=2^a3^b+1$, with $b\in\{1,2,3\}$, be an multiplicatively independent prime factor of $N_n$. Then
\begin{equation}\label{eq:lambda-relations}
2^\lambda\equiv\varepsilon\pmod q,\qquad
(\lambda,\varepsilon)=
\begin{cases}
(n-3a,1),&b=1,\\
(2n-3a,-1),&b=2,\\
(n-a,1),&b=3.
\end{cases}
\end{equation}
Moreover, $0<|\lambda|\le2n$.
\end{lemma}

\begin{proof}
The relations $q=2^a3^b+1$ and $N_n=2^n3^3+1$, together with $q\bigl\vert  N_n$, give $2^a3^b\equiv-1\pmod q$ and $2^n3^3\equiv-1\pmod q$. Since $N_n$ is odd and $N_n\equiv1\pmod3$, no prime divisor $q$ of $N_n$ can be $2$ or $3$. Thus $\gcd(q,6)=1$, so both $2$ and $3$ are invertible modulo $q$.

If $b=1$, cubing the first congruence and substituting into the second gives
$2^{n-3a}\equiv1\pmod q$. If $b=2$, cubing the first congruence and squaring the second gives $2^{2n-3a}\equiv-1\pmod q$. If $b=3$, direct substitution gives $2^{n-a}\equiv1\pmod q$. These are exactly the three cases in \eqref{eq:lambda-relations}.


Regardless of the choice of $b$, if $\lambda = 0$, then \eqref{eq:lambda-relations} shows that $3a = nb$, which shows that these terms are multiplicatively dependent \eqref{eq:dependence-condition}.
Hence, $\lambda\ne0$. Since $1\le a\le n$, the three cases also give $|\lambda|\le2n$.
\end{proof}


Squaring \eqref{eq:lambda-relations} gives $2^{2\lambda}\equiv1\pmod q$, and implies
\begin{equation}\label{eq:order-bound}
\ord_q(2)\bigl\vert 2|\lambda|,\qquad \ord_q(2)\le2|\lambda|\le4n.
\end{equation}
On the other hand, we have $\ord_q(2)\bigl\vert q-1$, while Korselt's criterion implies $q-1\bigl\vert 2^n3^3$. Putting together, we see that 
\begin{equation}\label{eq:order-form}
\ord_q(2)=2^\gamma3^s,\qquad 0\leq \gamma\leq n\text{ and }s\in\{0,1,2,3\}.
\end{equation}

With these relations at hand, we derive a general bound for divisors $M\bigl\vert 2^nk+1$, with $k\ge3$ odd, satisfying $2^u\equiv1\pmod M$ for some small $u$. The argument refines the box-principle method of \cite[Lemma~4]{CLP}, also used in \cite[Section~2.1]{Luca1729}.

\begin{lemma}\label{lem:divisor-bound}
Let $k\ge3$ be odd, put $L:=\log_2k$, let $M$ be a positive divisor of $2^nk+1$, and let $u\ge1$ satisfy $2^u\equiv1\pmod M$. Put $B:=\lfloor\sqrt{u/L}\rfloor$ and $A:=\lfloor u/(B+1)\rfloor$. Then
\begin{equation}\label{eq:divisor-bound}
A+BL<2\sqrt{uL},\qquad M\le2^Ak^B+1.
\end{equation}
In particular, $M<2^{2\sqrt{uL}}+1$.
\end{lemma}

\begin{proof}
If $M=1$, there is nothing to prove, so assume $M>1$. Since $B+1>\sqrt{u/L}$, we have $A\le u/(B+1)<\sqrt{uL}$; also $B\le\sqrt{u/L}$ gives $BL\le\sqrt{uL}$. Hence, $A+BL<2\sqrt{uL}$. Moreover, $A+1>u/(B+1)$, so $(A+1)(B+1)>u$.

Since $M\bigl\vert 2^nk+1$, we have $\gcd(M,2k)=1$ and $k\equiv-2^{-n}\pmod M$. Thus, for $(\alpha,\beta)\in\{0,\dots,A\}\times\{0,\dots,B\}$,
$$
2^\alpha k^\beta\equiv(-1)^\beta2^{\alpha-n\beta}\pmod M.
$$
There are $(A+1)(B+1)>u$ such pairs but only $u$ residue classes modulo $u$. Hence, by the pigeonhole principle, there exist two distinct pairs $(\alpha_1,\beta_1)$ and $(\alpha_2,\beta_2)$ that satisfy $\alpha_1-n\beta_1\equiv\alpha_2-n\beta_2\pmod u$. Since $2^u\equiv1\pmod M$, it follows that $(-1)^{\beta_1}2^{\alpha_1}k^{\beta_1}\equiv(-1)^{\beta_2}2^{\alpha_2}k^{\beta_2}\pmod M$.

After cancelling the common powers of $2$ and $k$, which are invertible modulo $M$, we obtain $M\bigl\vert  W$, where
$$
W:=2^{x_1}k^{y_1}-\delta2^{x_2}k^{y_2},\qquad \delta:=(-1)^{\beta_1+\beta_2}\in\{\pm1\},
$$
with $0\le x_i\le A$, $0\le y_i\le B$, $\min(x_1,x_2)=0$ and $\min(y_1,y_2)=0$.

We have $W\ne0$. Indeed, this is immediate if $\delta=-1$. If $\delta=1$ and $W=0$, then $2^{x_1}k^{y_1}=2^{x_2}k^{y_2}$; since $\gcd(2,k)=1$ and $k\ge3$, this forces $x_1=x_2$ and $y_1=y_2$, contradicting the choice of distinct pairs.

While bounding $|W|$, we encounter two cases up to switching the order of the pairs $(x_1,y_1)$ and $(x_2,y_2)$.
If $x_1=y_1=0$, then $|W|\le1+2^Ak^B$. Otherwise, we may assume $x_1=0$ and $y_2=0$, and then $|W|\le k^B+2^A\le1+2^Ak^B$, since $(2^A-1)(k^B-1)\ge0$. Thus $0<|W|\le1+2^Ak^B$. Since $M\bigl\vert  W$, we obtain $M\le2^Ak^B+1$, and \eqref{eq:divisor-bound} gives $M<2^{2\sqrt{uL}}+1$.
\end{proof}

For $k=27=3^3$, the following is immediate.

\begin{corollary}\label{cor:independent-factor-bound}
If $q=2^a3^b+1$ is an independent prime factor of the Carmichael number $N_n=27\,2^n+1$, then
$$
a+b\log_23<4\sqrt{n\log_227}.
$$
In particular, $a<4\sqrt{n\log_227}$.
\end{corollary}

\begin{proof}
Apply \autoref{lem:divisor-bound} with $k=27$, $M=q$ and $u=\ord_q(2)$. Since \eqref{eq:order-bound} gives $u\le4n$, we have $q-1\le2^A27^B<2^{2\sqrt{u\log_227}}\le2^{4\sqrt{n\log_227}}$. Taking base 2 logarithms yields the stated inequality.
\end{proof}

We now bound the product $R$ of all independent prime factors of $N_n$.

\begin{lemma}\label{lem:R-bound}
For each $s\in\{0,1,2,3\}$ with $3^s\le4n$, let $u_s$ be the largest integer of the form $2^\gamma3^s$ not exceeding $4n$, and set
\begin{equation}\label{eq:G-s}
G_s:=\gcd(2^{u_s}-1,N_n).
\end{equation}
If $3^s>4n$, set $G_s:=1$. Then every independent prime factor of $N_n$ divides at least one of $G_0,G_1,G_2,G_3$, and
\begin{equation}\label{eq:R-bound}
R\le\prod_{s=0}^3G_s,\qquad R<2^{16\sqrt{n\log_227}+4}.
\end{equation}
\end{lemma}

\begin{proof}
Let $q$ be an independent prime factor of $N_n$. By \eqref{eq:order-form}, $\ord_q(2)=2^\gamma3^s$ for some $0\leq \gamma\leq n$ and $s\in\{0,1,2,3\}$, and \eqref{eq:order-bound} gives $\ord_q(2)\le4n$. Hence, $3^s\le4n$, so $u_s$ is defined. By maximality of the exponent of $2$ in $u_s$, we have $\ord_q(2)\bigl\vert  u_s$. Therefore $q\bigl\vert 2^{u_s}-1$, and since $q\bigl\vert  N_n$, we obtain $q\bigl\vert  G_s$.

Thus every prime factor of $R$ divides at least one $G_s$. Since $R$ is squarefree, $R\bigl\vert \prod_{s=0}^3G_s$, proving the first inequality in \eqref{eq:R-bound}. For every $G_s>1$, \autoref{lem:divisor-bound} with $M=G_s$, $k=27$ and $u=u_s\le4n$ gives $G_s<2^{4\sqrt{n\log_227}}+1<2^{4\sqrt{n\log_227}+1}$. There are at most four such factors, which gives the second inequality in \eqref{eq:R-bound}.
\end{proof}

\begin{proposition}\label{prop:n-bound}
If $N_n=27\,2^n+1$ is Carmichael, then $n\le2742$. If, in addition, $3\nmid n$, then $n<1215$.
\end{proposition}

\begin{proof}
Let $y:=\sqrt n$ and $c:=16\sqrt{\log_227}$.

If $3\nmid n$, then \eqref{eq:R-versus-N} and \autoref{lem:R-bound} give $N_n-1=27\,2^{y^2}<N_n=R<2^{cy+4}$. Taking logarithms to base $2$ gives $y^2-cy+\log_227-4<0$. 
The larger root of the corresponding quadratic is approximately $34.8675$. The inequality and the larger root imply $n=y^2<1215.741$. Hence, $n\le1215$; since $3\nmid n$, the value $1215$ is excluded and therefore $n<1215$.

If $3\bigl\vert  n$, then \eqref{eq:R-versus-N} gives $R\ge N_n/D_n$. Put $x:=\sqrt[3]{N_n-1}=3\,2^{n/3}$. Since $N_n=x^3+1=(x+1)(x^2-x+1)$ and $D_n=x+1$, we have $R\ge x^2-x+1>x^2/2=9\,2^{2n/3-1}$. Combining this with \autoref{lem:R-bound} gives
$$
9\,2^{2y^2/3-1}<2^{cy+4}.
$$
Taking logarithms to base $2$ yields $(2/3)y^2-cy+\log_{2}{9}-5<0$. Similarly, the larger root is approximately $52.3861$, and we get $n=y^2<2744.303$. Since $3\bigl\vert  n$, it follows that $n\le2742$.
\end{proof}

Since $n$ is now bounded, \autoref{cor:independent-factor-bound} gives an explicit bound on the exponent $a$ of every independent prime factor. If $3\nmid n$, then $n<1215$, hence $a<4\sqrt{1215\log_2 27}<304.04$, so $a\le303$. If $3\bigl\vert  n$, then $n\le2742$, hence $a<4\sqrt{2742\log_2 27}<456.74$, so $a\le456$. In the latter case, there may also be the unique dependent prime factor $D_n=3\, 2^{n/3}+1$. Notice that any multiplicatively dependent term candidate $D_n$ is known explicitly and also bounded since $n$ is bounded.

\subsection{Comparison with Earlier Bounds}\label{rem:bounds}
Here we would like to highlight how the progress on the bounds evolved.

The general bound of \cite[Theorem~1]{CLP},
$$
n<2^{\,2\,10^7\tau(k)^2(\log k)^2\omega(k)},
$$
is effective but very large in the present case. For $k=27$, where $\tau(27)=4$ and $\omega(27)=1$, it gives $n<2^{\,2\,10^7\,16(\log27)^2}<2^{3.48\,10^9}$.

A much smaller bound can be extracted from the box principle argument of \cite[Section~2.1]{Luca1729}, which adapts \cite[Lemma~4]{CLP}. After separating the possible dependent factor, the estimate for the independent part applies equally to $p=3$. For $n\ge7$ the argument gives, for each of the four order classes,
$$
G_s\le2^{6\sqrt{n\log3}}e^{3\sqrt{n\log3}}+1<2^{10.5\sqrt{n\log3}},
$$
and hence $R<2^{42\sqrt{n\log3}}$. If there is no dependent factor, then $R=N_n$, so $n+\log_227<42\sqrt{n\log3}$, which gives $n\le1928$.

If a dependent factor occurs, then $3\bigl\vert  n$ and $D_n=3\,2^{n/3}+1\le(7/2)2^{n/3}$. In this case $N_n=D_nR$, so $\frac{2n}{3}+\log_2(54/7)<42\sqrt{n\log3}$, which gives $n\le4351$ and hence, since $3\bigl\vert  n$, $n\le4350$. The values $n\le6$ are already contained in this range.

Similar arguments can lead to a fixed bound for the exponent of an independent prime factor. By \cite[Lemma~4]{CLP}, if $q=2^ad+1$ divides a Carmichael number $2^nk+1$, with $d\bigl\vert  k$, $d>1$, and $2^ad$ multiplicatively independent of $2^nk$, then $a<7\sqrt{n\log k}$ whenever $n>3\log k$. For $k=27$ and $n\le4350$, this gives $a<7\sqrt{4350\log27}<838.16$, hence $a\le838$. If $n\le3\log27$, then $n\le9$ and trivially $a\le n\le9$. Thus $a\le838$ covers all independent prime factors in the range $n\le4350$.

\autoref{lem:divisor-bound} may be viewed as the same box principle argument with the two sides of the box chosen asymmetrically and optimised with respect to the relative sizes of $2$ and $k$. For $k=27$, since $\log_227=3\log3/\log2$, we have $4\sqrt{n\log_227}<8.33\sqrt{n\log3}$. Thus the new argument gives $G_s<2^{8.33\sqrt{n\log3}}+1$, improving the estimate on $G_s$ above. Likewise, \autoref{cor:independent-factor-bound} gives $a<4\sqrt{n\log_227}<4.81\sqrt{n\log27}$ instead of $a<7\sqrt{n\log27}$. Together with the treatment of the dependent factor, these improvements reduce the global range to $n\le2742$ and give $a\le456$ for every independent prime factor.

For comparison, the explicit estimate $n<25m^4\log p$ in \cite[Theorem~1.1(ii)]{Tiebekabe} is obtained under the hypothesis $p>m$. It therefore does not cover the exceptional cubic case considered here, where $p=m=3$.

\section{A Factorisation Free Sieve}\label{sec:sieve}

For each admissible $n$, define $a_{\max}(n):=\min\{n,\left\lceil4\sqrt{n\log_227}\, \right\rceil-1\}$ and
$$
\mathcal C_n:=
\left\{
2^a3^b+1:
\begin{array}{l}
1\le a\le a_{\max}(n),\ 1\le b\le3,\\
2^a3^b+1<N_n,\ 2^a3^b+1\bigl\vert  N_n
\end{array}
\right\}.
$$
If $3\bigl\vert  n$, we additionally include $D_n=3\,2^{n/3}+1$ in $\mathcal C_n$.

\begin{proposition}\label{prop:sieve}
If $N_n$ is Carmichael, then every prime divisor of $N_n$ belongs to $\mathcal C_n$. Furthermore, if $P_n:=\prod_{d\in\mathcal C_n}d$, then $N_n\bigl\vert  P_n$.
\end{proposition}

\begin{proof}
Let $q\bigl\vert  N_n$ be prime. By \eqref{eq:factor} and \autoref{thm:no-fermat}, $q=2^a3^b+1$ with $1\le b\le3$.

Suppose first that $q-1$ and $N_n-1$ are multiplicatively independent. By \autoref{cor:independent-factor-bound}, $a<4\sqrt{n\log_227}$; since $a$ is an integer, $a\leq\left\lceil4\sqrt{n\log_227}\,\right\rceil-1$. As also $a\le n$, we have $a\le a_{\max}(n)$. Finally, $q\bigl\vert  N_n$ and $q<N_n$ because $N_n$ is composite, so $q\in\mathcal C_n$.

If $q-1$ and $N_n-1$ are dependent, then \autoref{lem:dependent} gives $3\bigl\vert  n$ and $q=D_n$, which is included in $\mathcal C_n$. Thus every prime divisor of $N_n$ belongs to $\mathcal C_n$.

Since $N_n$ is Carmichael, it is square-free. Hence, $N_n$ is the product of its distinct prime divisors. Every such prime belongs to $\mathcal C_n$ and therefore divides $P_n$, so $N_n\bigl\vert  P_n$.
\end{proof}

\begin{proof}[Proof of \autoref{thm:taxicab}]
\autoref{prop:sieve} gives a finite sieve requiring neither factorisation of $N_n$ nor primality tests for the elements of $\mathcal C_n$: any admissible $n$ for which $N_n\nmid P_n$ cannot yield a Carmichael number. By \autoref{prop:n-bound}, it suffices to consider $1\le n<1215$ with $3\nmid n$, and $1\le n\le2742$ with $3\bigl\vert  n$. For each such $n$ we computed $P_n\bmod N_n$ with exact integer arithmetic. The only surviving value is
$$
n=6,
$$
for which
$$
\mathcal C_6=\{7,13,19\},\qquad N_6=27\,2^6+1=7\,13\,19=1729.
$$
As verified above by Korselt's criterion, $1729$ is Carmichael. This completes the proof.
\end{proof}

As an independent consistency check, we also ran the sieve with the earlier bounds of \autoref{sec:bound}, over $1\le n\le4350$ with $1\le a\le n$ and with $D_n$ included whenever $3\bigl\vert  n$. Again, the only surviving value is $n=6$. The sharper bounds of \autoref{prop:n-bound} therefore reduce the computational search substantially without changing the final conclusion.

Finally, we prove \autoref{cor:cubic}.

\begin{proof}[Proof of \autoref{cor:cubic}]
If $2^np^3+1$ is Carmichael, then \autoref{thm:prime-part} gives $p=3$ and \autoref{thm:taxicab} gives $n=6$. Conversely, $1729=7\,13\,19$ is square-free and $6$, $12$ and $18$ all divide $1728$, so $1729$ is Carmichael by Korselt's criterion.
\end{proof}

\section{Outlook}

In \cite{Luca1729}, Luca studies primes $p>3$. 
In Section 2, He describes how $p=3$ leads to the emergence of multiplicatively dependent factors. However, it is known that there would be at most one dependent prime factor \cite{CLP}. 

In the taxicab number $1729=7\,13\,19$, the factor $13=D_6=3\,2^2+1$ is precisely the dependent prime factor (of \autoref{lem:dependent} since $1728=12^3$), while $7=2\,3+1$ and $19=2\,3^2+1$ are independent. We observe that the exact $x^3+1$ structure that we see in \autoref{sec2} occurs in another similar Carmichael number $46657=2^6\,3^6+1=36^3+1=13\,37\,97$. This was also noticed in \cite{Tiebekabe}. Here $37=36+1$ is the unique dependent factor. It is not clear if 46657 is the only one with the shape $2^a 3^6 +1$ in its family. However, it lies outside the odd-$m$ hypothesis of \cite{Luca2026} and \cite{Tiebekabe}, where they study number of the form $2^n p^m +1$ with prime $p$. We plan to investigate the even-$m$ cases more to possibly put tractable bounds on prime divisors of Carmichael numbers of the form $2^k 3^m +1$, and provide some more underlying structure to Carmichael numbers in the near future.


\bibliographystyle{plain}
\bibliography{references}

@article{Korselt,
  author  = {Korselt, Alwin Reinhold},
  title   = {Probl{\`e}me chinois},
  journal = {L'Interm{\'e}diaire des Math{\'e}maticiens},
  volume  = {6},
  pages   = {142--143},
  year    = {1899}
}

@article{CLP,
  author  = {Cilleruelo, Javier and Luca, Florian and Pizarro-Madariaga, Amalia},
  title   = {Carmichael numbers in the sequence $(2^{n}k+1)_{n\ge 1}$},
  journal = {Mathematics of Computation},
  volume  = {85},
  number  = {297},
  pages   = {357--377},
  year    = {2016}
}

@article{AL,
  author  = {Alahmadi, Adel and Luca, Florian},
  title   = {There are no {C}armichael numbers of the form $2^{n}p+1$ with $p$ prime},
  journal = {Comptes Rendus Math{\'e}matique},
  volume  = {360},
  pages   = {1177--1181},
  year    = {2022},
  doi     = {10.5802/crmath.393}
}

@article{LR,
  author  = {Luca, Florian and Randrianantenaina, Jean Lucien},
  title   = {There is no {C}armichael number of the form $2^{n}p^{2}+1$ with $p$ prime},
  journal = {Integers},
  volume  = {23},
  pages   = {A51},
  year    = {2023},
  note    = {15 pp.},
  doi     = {10.5281/zenodo.8174516}
}

@article{Luca1729,
  author  = {Luca, Florian},
  title   = {Yet another curious property of $1729$},
  journal = {The Ramanujan Journal},
  volume  = {70},
  number  = {4},
  pages   = {74},
  year    = {2026},
  doi     = {10.1007/s11139-026-01450-1}
}

@misc{Luca2026,
  author        = {Luca, Florian},
  title         = {On {C}armichael numbers of the form $2^{n}p^{m}+1$},
  year          = {2026},
  eprint        = {2606.28829},
  archivePrefix = {arXiv},
  primaryClass  = {math.NT},
  note          = {Preprint}
}

@misc{Tiebekabe,
  author        = {Tiebekabe, Pagdame},
  title         = {Carmichael numbers of the form $2^{n}p^{m}+1$ with odd exponent $m$: refinements of a theorem of {F}lorian {L}uca, a general independence lemma, and a computational study},
  year          = {2026},
  eprint        = {2609.25104},
  archivePrefix = {arXiv},
  primaryClass  = {math.NT},
  note          = {Preprint}
}
\end{document}